\documentclass[11pt,reqno]{amsart}
\usepackage{amsmath, amssymb, amsthm}
\usepackage{url}
\usepackage{mathrsfs}
\usepackage[ansinew]{inputenc}
\usepackage[breaklinks]{hyperref}
\usepackage{fancyhdr}

\usepackage{caption}

\usepackage{multirow}

\allowdisplaybreaks

\renewcommand{\(}{\left\(}
\renewcommand{\)}{\right\)}
\renewcommand{\[}{\left\[}
\renewcommand{\]}{\right\]}
\newtheorem{remark}[]{Remark}
\numberwithin{equation}{section}
 \theoremstyle{plain}
\newtheorem{theorem}{Theorem}[section]

\newcommand{\sm}{\left(\begin{smallmatrix}}
\newcommand{\esm}{\end{smallmatrix}\right)}

\newtheorem{corollary}[theorem]{Corollary}

\newtheorem{example}[]{Example}

   \makeatletter
\def\proof{\@ifnextchar[{\@oproof}{\@nproof}}
\def\@oproof[#1][#2]{\trivlist\item[\hskip\labelsep\textit{#2 \textbf{Proof of}\
#1.}~]\ignorespaces}
\def\@nproof{\trivlist\item[\hskip\labelsep\textit{Proof.}~]\ignorespaces}
\makeatother

\usepackage{color}
\usepackage{amsmath}

\definecolor{blue}{rgb}{0,0,1}
\definecolor{red}{rgb}{1,0,0}
\definecolor{green}{rgb}{0,.6,.2}
\definecolor{purple}{rgb}{1,0,1}

\long\def\red#1\endred{{\color{red}#1}}
\long\def\blue#1\endblue{{\color{blue}#1}}
\long\def\purple#1\endpurple{{\color{purple}#1}}
\long\def\green#1\endgreen{{\color{green}#1}}

\begin{document}
\title[Parity, differences and smallest part in partition theory]{On Partition Classes Arising from Parity, Differences, and Repeated Smallest Parts}

\author{Rahul Kumar}
\address{Department of Mathematics, Indian Institute of Technology, Roorkee-247667, Uttarakhand, India}
\email{rahul.kumar@ma.iitr.ac.in} 

\author{Nargish Punia}
\address{Department of Mathematics, Indian Institute of Technology, Roorkee-247667, Uttarakhand, India}
\email{nargish@ma.iitr.ac.in}

\maketitle
\pagenumbering{arabic}
\pagestyle{headings}
\begin{abstract}
In this paper, we study various classes of partition functions such as those related to the parity of the number of parts, to differences of partition numbers, and to partitions with a repeated smallest part. We establish identities connecting these various classes of partitions. Moreover, our identities help us to extend the Euler's partition theorem. An analogue of Legendre's theorem of the partition-theoretic interpretation of Euler's pentagonal number theorem is also derived. Both combinatorial and $q$-series proofs are given for our results.

\end{abstract}
%\tableofcontents

\section{Introduction}\label{intro}

Let $A(n)$ denote the number of partitions of $n$ into distinct parts, and let $B(n)$ denote the number of partitions of $n$ into odd parts. Then one of the most celebrated theorems in the theory of partitions is Euler's theorem, which states that
\begin{align*}
A(n)=B(n).
\end{align*}

Very recently, Andrews, Yee  and the first author \cite{aky} introduced two further important classes of partitions to Euler's theorem. More explicitly, they derived the following theorem.
\begin{theorem}\textup{(}\cite[Theorem 1]{aky}\textup{)}\label{main theorem}
For $n>0$, we have
\begin{align}\label{main theorem eqn}
A(n)=B(n)=C(n+1)=\frac{1}{2}D(n+1),
\end{align}
where $C(n)$ is the number of partitions of $n$ with largest part even and parts not exceeding half of the largest part are distinct, and $D(n)$ is the number of partitions of $n$ into non-negative parts wherein the smallest part appear exactly twice and no other parts are repeated.
\end{theorem}

Analogues of this theorem in two different directions related to Glaisher's partition theorem \cite{glaisher} are subsequently obtained by Andrews and Dhar \cite{ad} and by Lin and Zang \cite{lz}. 

Note that the partition function $D(n)$ counts the number of partitions of $n$ into non-negative parts wherein the smallest part appear exactly twice and no other parts are repeated. Now this naturally raises the question that what happens if the smallest part is allowed to appear an arbitrary number of times, say $k$, while all other parts remain distinct? Let us denote the number of such partitions by  $D_k(n)$.  This question serves as the starting point of the present paper, and it will soon become evident that it opens up several further directions to explore.

Before stating our results, we introduce the following partition functions based on the parity. 
%Let $C_{k}^e(n)$ \textup{(}resp. $C_{k}^o(n)$\textup{)}  denote the number of partitions of $n$  satisfying the conditions:
%
%(1) the partition contains an even part \textup{(}say $2\ell$\textup{)} and all parts not exceeding $\ell$ are distinct; and 
%
%(2) there are at most $(k-1)$ distinct even parts in the interval $[2\ell+2,2\ell+2k-2]$ such that the number of such parts is even \textup{(}resp. odd\textup{)}.
Let  $B_{k}^e(n)$ \textup{(}resp. $B_{k}^o(n)$\textup{)} denote the number of partitions of $n$ in which all parts are odd, except possibly for at most $k-1$ distinct even parts lying in the interval $[2\ell+2,2\ell+2k-2]$, where $2\ell-1$ is the largest odd part, the number of such distinct even parts being even (resp. odd).

Moreover, let $C_{k}^e(n)$ \textup{(}resp. $C_{k}^o(n)$\textup{)}  denote the number of partitions of $n$ in which the partition contains an even part (say $2\ell$), all parts not exceeding $\ell$ are distinct, and the only parts larger than $2\ell$ are at most $k-1$  distinct even parts lying in $[2\ell+2,2\ell+2k-2]$, the number of such distinct even parts being even (resp. odd).

% Moreover, let  $B_{k}^e(n)$ \textup{(}resp. $B_{k}^o(n)$\textup{)}  denote  the number of partitions of $n$ satisfying the conditions:
% 
% (1)  the odd parts are unrestricted (but say the largest part is $\ell$); and 
% 
% (2) there are at most $(k-1)$ distinct even parts lying  in the interval $[\ell+3,\ell+2k-1]$, and the number of such parts is even (resp. odd). 

Our first result is stated below.
\begin{theorem}\label{B_k=C_k}
For any positive integers $n\geq1$ and $k\geq1$, we have
\begin{align}
B_{k}^e(n)=C_{k}^e(n+1),\nonumber
\end{align}
and 
\begin{align}
B_{k}^o(n)=C_{k}^o(n+1).\nonumber
\end{align}
\end{theorem}

%\begin{remark}
%Observe that the structure of $B_k^o(n)$ is quite different from $B_k^e(n)$ in the sense that there would always be atleast one part  coming from the interval $[2\ell+2,2\ell+2k-2]$ in $B_k^o(n)$ as number of parts is required to be odd and therefore cannot be zero, whereas this is not the case with $B_k^e(n)$. The same goes for $C_k^o(n)$ and $C_k^e(n)$.
%\end{remark}

\begin{example}
We take $n=8$ and $k=2$ to illustrate Theorem \ref{B_k=C_k}.
\begin{center}
\renewcommand{\arraystretch}{1.3}
\begin{tabular}{c @{\hspace{1cm}} c @{\hspace{1cm}} c @{\hspace{1cm}} c}
\textbf{$B_{2}^e(8)$} & \textbf{$C_{2}^e(9)$} & \textbf{$B_{2}^o(8)$} & \textbf{$C_{2}^o(9)$} \\[2mm]
\hline \\[-2mm]
$7+1$ & $8+1$ & $4+1+1+1+1$ & $4+2+2+1$ \\
$5+3$ & $6+3$ &  &  \\
$5+1+1+1$ & $6+2+1$ &  &  \\
$3+3+1+1$ & $4+4+1$ &  &  \\
$1+1+1+1+1+1+1+1$ & $4+3+2$ &  &  \\
$3+1+1+1+1+1$ & $2+2+2+2+1$ &  & 
\end{tabular}
\end{center}
	\medskip
	
This table shows that  $B_2^e(8)=6=C_2^e(9)$ and $B_2^o(8)=1=C_2^o(9)$.
\end{example}

Differences of partition numbers have been studied extensively in the literature, see, for example, \cite{andrews difference,  ay, bach}, among others. In particular, differences of partition numbers according to parity are of fundamental importance (see \cite{andrews ac} and references therein). A classical example, due to Legendre \cite[pp.~128--133]{legendre} (\cite[p.~10, Theorem 1.6]{andrews}), arises from an interpretation of Euler's pentagonal number theorem, namely,
\begin{align}\label{legendre}
P_e(d,n)-P_o(d,n)=
\begin{cases}
(-1)^m, & \textup{if}\ n=\frac{1}{2}m(3m\pm1), \\
0, & \text{otherwise},
\end{cases}
\end{align}
where $P_e(d,n)$ (resp. $P_o(d,n)$) denote the number of partitions
of $n$ into an even (resp. odd) number of distinct parts.

In this direction, our next result shows that differences of the partition numbers appearing in Theorem \ref{B_k=C_k} are precisely `equal' to $D_k(n)$, the number of partitions into non-negative parts wherein smallest part is allowed to repeat exactly  $k$ times, while all other parts remain distinct, as defined earlier.
\begin{theorem}\label{BkCkDk}
Let $k$ be any natural number. For all $n\geq2^{k-1}k(2k-1)$, we have
\begin{align}\label{BkCkDk eqn}
B_{k}^e(n)-B_{k}^o(n)=C_{k}^e(n+1)-C_{k}^o(n+1)=\frac{1}{2}D_{2k}(n+1).
\end{align}
\end{theorem}

\begin{remark}
Equation \eqref{BkCkDk eqn} immediately  implies that $B_{k}^e(n)>B_{k}^o(n)$ and $C_{k}^e(n)>C_{k}^o(n)$ for $n\geq2^{k-1}k(2k-1)$,  a fact that does not seem to be easy to establish directly.
\end{remark}

Our next result connects $D_k(n)$ with another type of partitions.
\begin{theorem}\label{AkDk}
Let $k$ and $n$ be any natural numbers. Then, we have
\begin{align}\label{AkDk eqn}
A_k(n)=\frac{1}{2}D_{k}(n+1),
\end{align}
where $
A_k(n):=\frac{1}{2}\left(P_1(d,n)+P_{k-1}'(d,n)+P_2(d,n)+P_{k-1}^{''}(d,n)\right).$
Here, and throughout the paper, $P_1(d,n)$ counts partitions of $n$ into distinct parts whose smallest part exceeds 
$1$ and $P_2(d,n)$ denote the number of partitions of $n$ into distinct parts for which the difference between the two smallest parts is at least $2$. Moreover the function $P_{k-1}'(d,n)$  enumerates partitions of  $n$ in which all parts are distinct except that the smallest part must be $1$ with multiplicity $k-1$; and finally $P_{k-1}^{''}(d,n)$ denote the number of partitions of $n$ into distinct parts  except that the second smallest part appears with multiplicity $k-1$, and the difference between the two smallest parts is $1$.				
\end{theorem}

Our next result unifies Theorem \ref{BkCkDk} and Theorem \ref{AkDk} to produce a generalization of Theorem \ref{main theorem} of Andrews, Kumar, and Yee \cite{aky}.
\begin{theorem}\label{genral result}
Let $k$ be any natural number. For all $n\geq2^{k-1}k(2k-1)$, we have
\begin{align}\label{genral result eqn}
A_{2k}(n)=B_{k}^e(n)-B_{k}^o(n)=C_{k}^e(n+1)-C_{k}^o(n+1)=\frac{1}{2}D_{2k}(n+1).
\end{align}
\end{theorem}
			
The special case $k=1$ of Theorem \ref{genral result} reduces to Theorem \ref{main theorem}, as shown in the following corollary.
\begin{corollary}\label{special case k=1}
The equation \eqref{main theorem eqn} holds true.
\end{corollary}

Our next result extends Theorem~\ref{main theorem} by incorporating two additional partition functions.
\begin{theorem}\label{extd euler}
For $n>0$, we have
\begin{align}
A(n)=B(n)=C(n+1)=\frac{1}{2}D(n+1)=E(n+2)=F(n+1),\nonumber
\end{align}
where $E(n)$ denotes the number of partitions of $n$ having all parts odd with unique largest part and $F(n)$ denotes the number of partitions of $n$ with unique largest even part and all other parts are odd.			 
\end{theorem}

Continuing our discussion of differences of partition numbers and parity of parts, we now present an elegant analogue of Legendre's assertion \eqref{legendre} in the next result.
\begin{theorem}\label{legendre analogue}
Let $D_k^e(n)$ \textup{(}resp. $D_k^o(n)$\textup{)} denote the number of partitions counted by $D_k(n)$ with the additional condition that the number of parts greater than the smallest part is even \textup{(}resp. odd\textup{)}. Then
\begin{align}\label{legendre analogue eqn}
D_k^e(n)-D_k^o(n)=
\begin{cases}
(-1)^m, & \textup{if}\ 1\leq n=\frac{1}{2}m(3m\pm1)\leq k-1, \\
P_e^{(k)}(d,n)-P_o^{(k)}(d,n), & \textup{if}\ k-1<n\leq\frac{1}{2}k(k-1),\\
0, & \text{otherwise},
\end{cases}
\end{align}
where $P_e^{(k)}(d,n)$ \textup{(}resp. $P_o^{(k)}(d,n)$\textup{)} denote the number of partitions of $n$ into an even \textup{(}resp. odd\textup{)} number of distinct parts in which no part exceeds $k-1$.
\end{theorem}

The theorem above gives the following important corollary.
\begin{corollary}\label{parity of D_k}
For $n>\frac{k(k-1)}{2}$, among the partitions counted by $D_k(n)$, those with an even number of parts exceeding the smallest part are equinumerous with those having an odd number of such parts.

Consequently, $D_k(n)$ is always even for all  $n>\frac{k(k-1)}{2}$, where $k$ is any  positive integer.
\end{corollary}

\begin{remark}
Note that while Theorem \textup{\ref{main theorem}} shows that $D_2(n)$ \textup{(}equivalently, $D(n)$\textup{)} is even, Corollary \textup{\ref{parity of D_k}} gives another proof of this fact. In fact, it reveals a stronger structural result, namely that for each $k$, $D_k(n)$ $($and hence $D_2(n))$ decomposes naturally into two disjoint partition classes.
\end{remark}

\begin{example}
Consider the case $D_2(7)$ to demonstrate the above result. The admissible partitions of $7$ counted by $D_2(n)$ are 
$$0+0+7
,\ 0+0+6+1,\ 0+0+5+2,\ 0+0+4+3,\ 0+0+4+2+1,\ 1+1+5,\ 1+1+2+3,\  2+2+3.$$ 
Among these,  the four partitions 
$$0+0+6+1,\ 0+0+5+2,\ 0+0+4+3,\ 1+1+2+3$$ 
have even number of parts greater than the smallest part, whereas the remaining four partitions
  $$0+0+7,\ 0+0+4+2+1,\ 1+1+5,\ 2+2+3$$
 have an odd number of parts greater than the smallest part. Hence, both classes are equinumerous.
\end{example}

\begin{remark}
In the view of Corollary \textup{\ref{parity of D_k}} or Theorem \textup{\ref{legendre analogue}}, we can say that Theorem \textup{\ref{main theorem}} above can be reformulated in a more refined way, namely,
\begin{align}\label{main theorem re}
A(n)=B(n)=C(n+1)=\frac{1}{2}D(n+1)=D_2^e(n+1)=D_2^o(n+1).
\end{align}
More generally, our Theorem \textup{\ref{genral result eqn}} can be recast in the following way:
\begin{align}\label{genral result re}
A_{2k}(n)=B_{k}^e(n)-B_{k}^o(n)=C_{k}^e(n+1)-C_{k}^o(n+1)=\frac{1}{2}D_{2k}(n+1)=D_{2k}^e(n+1)=D_{2k}^o(n+1).
\end{align}
\end{remark}

Before stating our further results, we note that Andrews and Bachraoui \cite{ab} recently studied the function spt$k_d(n)$, which counts the number of partitions of $n$ in which the smallest part occurs exactly $k$ times and the remaining parts do not repeat. Note that our function $D_k(n)$ differs from spt$k_d(n)$ only in that $D_k(n)$ allows non-negative parts as well. We next derive some properties of $D_k(n)$. It is clear that its generating function can be given by\footnote{Here we use the convention that $D_k(0)=1$.}
\begin{align}\label{generating function of D_k}
\sum_{n=0}^{\infty}D_k(n)q^n=\sum_{n=0}^\infty q^{nk}(-q^{n+1};q)_\infty,
\end{align}
here, and throughout the sequel,
\begin{align*}
(A;q)_N=(1-A)(1-Aq)\cdots\left(1-Aq^{N-1}\right) \qquad (N\in\mathbb{N}),
\end{align*}
and
\begin{align*}
(A;q)_0=1,
\end{align*}
and 
\begin{align*}
(A;q)_\infty=\lim_{N\to\infty}(A;q)_N.
\end{align*}

We here note that the study of integer partitions with a repeating smallest part has received considerable attention in recent years; see, for example, \cite{ab, ad, bm, dks}, among others.

In the next theorem, we show that the generating function of $D_k(n)$ in \eqref{generating function of D_k} can be written explicitly as follows.
\begin{theorem}\label{D_k in terms of polyn}
For any positive integer $k$, we have
\begin{align}\label{D_k in terms of polyn eqn}
\sum_{n=0}^{\infty}D_k(n)q^n&=\sum_{n=0}^\infty q^{nk}(-q^{n+1};q)_\infty\nonumber\\
&=2(-q;q)_\infty\sum_{j=0}^{k-1}(-1)^j(q^{k-j};q)_j+(-1)^k(q;q)_{k-1}.
\end{align}
\end{theorem}

The rightmost identity in the above result follows directly from its finite analogue $q$-series identity, stated in the next theorem.
\begin{theorem}\label{finite analogue}
If $k\ge 1$ and $N\ge 0$ then
\begin{align}\label{finite analogue eqn}
\sum_{j=0}^{N}q^{kj}(-q^{j+1};q)_\infty
%=(-q;q)_\infty\left(2+2\sum_{j=1}^{k-1}(-1)^j(q^{k-j};q)_j+\frac{(-1)^k(q;q)_{k-1}}{(-q;q)_N}\right).
=(-q;q)_\infty\sum_{j=0}^{k-1}(-1)^{j+k-1}\left(q^{j+1};q\right)_{k-j-1}\left(2-\frac{q^{(N+1)j}}{(-q;q)_N}\right).
\end{align}
\end{theorem}
The special case $k=1$ of this identity is the following result of Andrews \cite[p.~919, (2.4)]{andrews amm}:
\begin{align}\label{andrews amm}
\sum_{j=0}^{N}\frac{q^{j}}{(-q;q)_j}=2-\frac{1}{(-q;q)_N}.
\end{align}

The result in Theorem \ref{D_k in terms of polyn} produces the following identity.
\begin{theorem}\label{D_k+D_k}
For $n>k-1$, we have
\begin{align}\label{D_k+D_k  eqn}
D_{k}(n)+D_{k-1}(n)=D_{k-1}(n-k+1)+2A(n).
\end{align}						
\end{theorem}

%If we take $k=2$ in Theorem \ref{D_k in terms of polyn} or Theorem \ref{D_k+D_k} then it gives one of the identities from Theorem \ref{main theorem}:
%\begin{corollary}\label{A_n=D_n}
%For any natural number $n$, we have
%\begin{align*}
%A(n)=\frac{1}{2}D(n+1),
%\end{align*}
%where $A(n)$ and $D(n)$ are as defined in Theorem \textup{\ref{main theorem}}.
%\end{corollary}

We conclude this section by showing that Theorem \ref{D_k in terms of polyn} produces infinitely many results of the form presented below.
\begin{corollary}\label{D3A}
For $n\geq4$, we have
\begin{align}\label{D3A eqn}
D_3(n)=2A(n-3)-2A(n-1)+2A(n).
\end{align}
\end{corollary}

This paper is organized as follows. Section \ref{bijective} is devoted to combinatorial proofs. The $q$-series proofs of all the results stated in Section \ref{intro} are given in Section \ref{analytic}. We conclude the paper with some remarks and further questions in Section \ref{concluding remarks}.
\section{Bijective Proofs}\label{bijective}

\subsection{A bijective proof of $A_k(n)=\frac{1}{2}D_k(n+1)$}
We first provide a bijection between the sets containing the partitions enumerated by $D_k(n)$ and $P_1(d,n-1)+P_2(d,n-1)+P_{k-1}'(d,n-1)+P_{k-1}''(d,n-1),\ n>1$.

Consider a partition $\lambda=(\lambda_1,\lambda_2,\dots,\lambda_\ell)$ counted by $D_{k}(n)$. Since $D_{k}(n)$ can be viewed as counting partitions of $n$ either into distinct parts or with the smallest part appearing exactly $k$ times and and other parts remains distinct, we accordingly distinguish between these two cases, namely, 
\begin{enumerate}
\item[] Case 1: all parts $\lambda_i,\ 1\leq i\leq\ell$, are distinct, or
\item[] Case 2: the smallest part is repeated $k$ times and other parts are distinct.
\end{enumerate}

\begin{enumerate}
\item[] \textbf{Case 1:} $ \lambda_1>\lambda_2>\dots>\lambda_\ell>0$. In this situation, we further distinguish two cases: $\lambda_\ell>1$ and $\lambda_\ell=1$. 

\begin{enumerate}
\item[] \textbf{Subcase 1:} First, suppose that $\lambda_\ell>1$. Subtracting $1$ from the smallest part $\lambda_\ell$ produces a partition of $n-1$ into $\ell$ distinct parts in which the difference between the two smallest parts is at least $2$. Hence, the resulting partition is counted by $P_2(d,n-1)$.

\item[] \textbf{Subcase 2:} Next consider the case  $\lambda_\ell=1$. We subtract 1 from the smallest part $\lambda_\ell$. In this process, we will get a partition of $n-1$ into $\ell-1$ distinct parts where the smallest part is greater than $1$. Thus, the resulting partition is enumerated by $P_1(d,n-1)$.
\end{enumerate}

\item[] \textbf{Case 2:} $	\lambda_1>\lambda_2>\dots>\lambda_{(\ell-(k-1))}=\lambda_{(\ell-(k-2))}=\dots=\lambda_\ell>0$. Similar to the previous case, we again distinguish two cases: $\lambda_\ell>1$ and $\lambda_\ell=1$.
\begin{enumerate}
\item[] \textbf{Subcase 1:} First, suppose that $\lambda_{(\ell-(k-1))}=\lambda_{(\ell-(k-2))}=\dots=\lambda_{\ell-1}=\lambda_\ell>1$. Upon subtracting 1 from $\lambda_\ell$, we get a partition  counted by $P_{k-1}''(d,n-1)$ as all the parts are distinct except that the second smallest part appears $k-1$ times and the difference between the two smallest parts is 1.
\item[] \textbf{Subcase 2:}
Now consider the case $\lambda_1>\lambda_2>\dots>\lambda_{(\ell-(k-1))}=\lambda_{(\ell-(k-2))}=\dots=\lambda_{\ell-1}=\lambda_\ell=1$. If we subtract $1$ from $\lambda_\ell$, then it leads to a partition counted by $P_{k-1}'(d,n-1)$ as all the parts are distinct except that the smallest part is $1$ with multiplicity $k-1$.
\end{enumerate}
\end{enumerate}
This completes the proof.

\begin{example}	
To illustrate the above bijection, we consider the case $n=8$ and $
k=4$. The table below shows how partitions counted by $D_4(8)$ are mapped to other classes of partitions under this bijection. The partition appearing next to each partition of $D_4(8)$ is its image. \\

%	\begin{tabular}{ccccc}\\
%		\textbf{$D_4(8)$} &\quad \textbf{$P_1(d,7)$}&\quad \textbf{$P_{k-1}(d,7)$}&\quad \textbf{$P_2(d,7)$}&\quad \textbf{$P_{k-1}^{'}(d,7)$}\\
%		8 &\quad 7 &\quad &\quad &\quad\\
%		6+2 &\quad 6+1&\quad&\quad &\quad\\
%		7+1 &\quad &\quad&\quad7\\
%		5+3 &\quad 5+2&\quad&\quad\\
%		5+2+1 &\quad&\quad &\quad 5+2\\
%		4+3+1 &\quad&\quad&\quad 4+3\\
%		4+1+1+1+1 &\quad&\quad 4+1+1+1&\quad\\
%		2+2+2+2 &\quad&\quad &\quad&\quad 2+2+2+1\\
%		\end{tabular}
%		
		\begin{center}
\renewcommand{\arraystretch}{1.2}
\begin{table}[htbp]
\centering
\renewcommand{\arraystretch}{1.2}
\begin{tabular}{c|cccc}
\hline
\textbf{$D_4(8)$} 
& \textbf{$P_1(d,7)$} \qquad\qquad
& \textbf{$P_{k-1}'(d,7)$} \qquad\qquad
& \textbf{$P_2(d,7)$} \qquad\qquad
& \textbf{$P_{k-1}''(d,7)$} \\
\hline
8              &            &          &\hspace{-1cm}7       &          \\
6+2            &            &          &\hspace{-1cm}6+1     &           \\
7+1            & \hspace{-1cm} 7          &          &         &           \\
5+3            &            &          &\hspace{-1cm}5+2     &           \\
5+2+1          & \hspace{-1cm}5+2        &          &         &           \\
4+3+1          &\hspace{-1cm} 4+3        &          &         &           \\
4+1+1+1+1      &            & 4+1+1+1   &        &           \\
2+2+2+2        &            &           &        & 2+2+2+1         \\
\hline
\end{tabular}
\end{table}
\end{center}

\end{example}

\medskip

\subsection{A bijective proof of $B_{k}^e(n)=C_{k}^e(n+1)$ and $B_{k}^o(n)=C_{k}^o(n+1)$}
Note that when $k=1$, $B_{1}^e(n)=B(n)$, $C_{1}^e(n+1)=C(n+1)$ and $B_{1}^o(n)=0$, $C_{1}^o(n)=0$, where $B(n)$ and $C(n)$ are as defined in Theorem \ref{main theorem}. A bijection between the sets enumerated by $B(n)$ and $C(n+1)$ was provided in \cite{aky} through Glaisher's bijection \cite{glaisher}. We will be referring to this bijection multiple times in our proof.

The proof is carried out by induction on $k$. The base case $k=1$ follows from the preceding discussion. We first describe the bijective map for the case $k=2$ before treating the general case.

\noindent
\textbf{The case $k=2$:} Observe that $B_{2}^e(n)=B(n)$ and $C_{2}^e(n+1)=C(n+1)$. Since Theorem \ref{main theorem} asserts that $B(n)=C(n+1)$, the bijection constructed in \cite{aky} provides a bijection between the sets of partitions counted by $B_{2}^e(n)$ and $C_{2}^e(n+1)$.

We next show a bijection between the sets of partitions enumerated by $B_{2}^o(n)$ and $C_{2}^o(n+1)$. To that end, let $\lambda$ be a partition counted by $C_{2}^o(n)$. Then it is of the form:
$$(2\ell+2)+(2\ell)\ +\ \textup{the parts not exceeding $\ell$ are distinct, while parts in  $(\ell,2\ell]$ are unrestricted}.$$
Now, clearly the partition,
$$(2\ell)+\ \textup{the parts not exceeding $\ell$ are distinct, while parts in  $(\ell,2\ell]$ are unrestricted}.$$
is counted by $C_{2}^e(n-2\ell-2)$.  Using the bijection between $B_{2}^e(n)$ and $C_{2}^e(n+1)$, this partition can be mapped to a partition $\mu$ counted by $B_{2}^e(n-2\ell-3)$. Appending the part $2\ell+2$ to $\mu$ then yields a partition  enumerated by $B_{2}^o(n-1)$. 
			
Conversely, consider a partition counted by $B_{2}^o(n)$. Such a partition is of the form
$$(2m)+(2m-3)+ \textup{ odd parts}.$$ 
Next observe that the partition,
$$(2m-3)+\textup{ odd parts}.$$ 
is counted by $B_{2}^e(n-2m)$. By applying the bijection between $B_{2}^e(n)$ and $C_{2}^e(n+1)$, this partition corresponds to a partition $\mu$ counted by $C_{2}^e(n-2m+1)$  . Finally, adjoining the part $2m$ to $\mu$ produces a partition counted by $C_{2}^o(n+1)$.

This provides the required bijection for the case $k=2$.

\noindent
\textbf{The general case $k$:}
Let $\lambda$ be a partition counted by $C_{k}^e(n+1)$, then $\lambda$ is of the form:
\begin{align*}
\underbrace{\textup{parts belonging to the set}\  S=\{2\ell+2k-2,2\ell+2k-4,\dots,2\ell+2\}}_{\text{number of parts are even}}\\ +\ 2\ell\ +\textup{the parts not exceeding}\ \ell\ \textup{are distinct, while parts in}\  (\ell,2\ell]\ \textup{are unrestricted}.
\end{align*}
First of all, note that if no part from the set $S$ appears in $\lambda$, then such a partition is precisely counted by $C(n+1)$. In this case, the known aforementioned bijection maps it to a partition counted by $B(n)$, and hence they are, in particular, counted by $B_k^e(n)$.

Now suppose that at least one part from $S$ appears in $\lambda$, and let $m\in S$ be the largest such part. Writing
			 $$m+ \textup{remaining parts.}$$
we see that the remaining parts contain an odd number of distinct parts from the set $S$. Consequently, the remaining parts form a partition counted by $C_{k-1}^o(n-m+1)$. By the induction hypothesis, this partition can be mapped to a partition $\mu$ counted by $B_{k-1}^o(n-m)$. Finally, adjoining the part $m$ to $\mu$ yields a partition counted by $B_k^e(n)$.

Conversely, let $\lambda$ be a partition counted by $B_{k}^e(n)$, then $\lambda$ is of the form,
\begin{align*}
&\underbrace{\text{parts belonging to the set}\ S=\{2l+2k-2,2l+2k-4,\dots,2l+2\}}_{\text{number of parts are even}}\\
&\qquad\qquad\quad+\ (2l-1)\ +\text{remaining parts are odd}.
\end{align*}
Similar to the previous case, if no part from the set $S$ appears in $\lambda$, then such a partition is precisely counted by $B(n)$. In this case, the known aforementioned bijection maps it to a partition counted by $C(n+1)$, and hence they are, in particular, counted by $C_k^e(n)$.

Now suppose that at least one part from $S$ appears in $\lambda$, and let $m\in S$ be the largest such part. Writing
			 $$m+ \textup{remaining parts}.$$
Clearly, these remaining parts represent a partition counted by $B_{k-1}^o(n-m)$. Now by the induction hypothesis, this partition can be mapped to a partition $\mu$ counted by $C_{k-1}^o(n-m+1)$. Finally, adjoining the part $m$ to $\mu$ yields a partition counted by $C_k^e(n+1)$.

In same way, one can give a bijective proof for $$B_{k}^o(n)=C_{k}^o(n+1).$$
This completes the proof.

\begin{example}	
For bijection between $B_{k}^e(n)$ and $C_{k}^e(n+1)$ for $k=3$. Consider a partition $(8,6,4,4)$ is a partition counted by $C_{3}^e(22)$. Clearly, in this case $2\ell=4$ and the subpartition partition $(6,4,4)$ is counted by $C_{2}^o(14)$. This subpartition can be mapped to the partition $(6,3,1,1,1,1)$ which is  counted by $B_{2}^o(13)$, by converting the pair $(4,4)$ into $(3,1,1,1,1)$ as follows:
\begin{itemize}
\item subtract 1 from the largest part $4$.
\item apply Glaisher's bijections to the remaining part $4$ by writing $4=2^2.1$ \textup{(}See \cite{aky, glaisher} for more details\textup{)}.
\end{itemize}
Finally, adjoining the removed part $8$ to this new partition, we obtain $(8,6,3,1,1,1,1)$ which is counted by $B_{3}^e(21)$.

Conversely, the partition $(8,6,3,1,1,1,1)$ which is counted by $B_{3}^e(21)$ can be mapped back to $(8,6,4,4)$ counted by $C_{3}^e(22)$. To that end, we convert the subpartition $(3,1,1,1,1)$ into $(4,4)$ using the bijection map between $B(n)$ and $C(n+1)$ as follows \cite{aky}:
\begin{itemize}
\item adding 1 to the largest part $3$, obtaining $4$,
\item here  $k(1)=2$ \textup{(}Note that $k(a)$ is the smallest integer such that $2^k>\frac{N}{a}$ and $\mu=2N-1$ is the largest odd part present in the partition and odd part $a\le N$ repeat with multiplicity $f$\textup{)}. In our case, the part $a=1$ occurs with multiplicity $f=4$ and the binary representation of 4 is $2^2$. Hence, the binary expansion yields \textup{(}for further details of the method see \cite{aky}\textup{)} 
					$$ (1,1,1,1)\longrightarrow4.$$
\end{itemize}
\end{example}	
				
\begin{example}	
Take $k=3$ and $n=7$. It is clear that $B_3^o(7)=1$, with the unique admissible partition being $4+1+1+1$, and that $C_3^o(8)=1$, whose only allowed partition is $4+2+2$. Hence, there is no alternative choice, these two partitions must correspond to each other under the bijection. This can also be verified directly using the bijection described above. 

For the even case, the bijection is described in the following table. The partition written in front of each partition denotes its image under the bijection.
\begin{table}[h]
\centering
\renewcommand{\arraystretch}{1.3}

\begin{tabular}{cc}

\textbf{$C_{3}^e(8)$} & \textbf{$B_{3}^e(7)$} \\
\hline
$8$ & $7$ \\

$6+2$ & $5+1+1$ \\

$4+4$ & $3+1+1+1+1$ \\

$2+2+2+2$ &\quad $1+1+1+1+1+1+1$ \\

\end{tabular}
\end{table}
	\end{example}
	
\subsection{\textbf{A bijective proof of Theorem \textup{\ref{D_k+D_k}}}}
Let $\lambda$ be a partition counted by $D_{k}(n)$ or counted by $ D_{k-1}(n)$.

\noindent
\textbf{Case 1:} if $\lambda$ is counted by  $D_{k}(n)$  then it is of the form $\lambda_1>\lambda_2>\dots>\lambda_{\ell-k+1}=\lambda_{\ell-k+2}=\dots=\lambda_{\ell+1}=\lambda_\ell$.
\begin{itemize}
\item If $\lambda_\ell=1$, we subtract $1$ from each of the first $k-1$ smallest parts to obtain a partition counted by $D_{k-1}(n-k+1)$, in which the smallest part is $0$, all parts are distinct, and the second smallest part is equal to $1$.
\item If $\lambda_\ell>1$, then subtracting $1$ from the first $k-1$ smallest parts yields a partition counted by $D_{k-1}(n-k+1)$, in which the two smallest parts are consecutive and the smallest part is at least $1$.
\item If $\lambda_\ell=0$ then clearly it is a partition counted by $A(n)$.
\end{itemize} 
\textbf{Case 2:} if $\lambda$ is counted by  $D_{k-1}(n)$, then it is of the form $\lambda_1>\lambda_2>\dots>\lambda_{\ell-k+2}=\dots=\lambda_{\ell+1}=\lambda_\ell$ (smallest part repeats $k-1$ times).
\begin{itemize}
\item If $\lambda_\ell=1$, we subtract $1$ from each of the first $k-1$ smallest parts to obtain a partition counted by $D_{k-1}(n-k+1)$,  where the smallest part is $0$ and all parts are distinct with second smallest part$>1$.
\item If $\lambda_\ell>1$, then subtracting $1$ each of the first $k-1$ smallest parts to obtain a partition counted by $D_{k-1}(n-k+1)$, whose two smallest parts are not consecutive.
\item If $\lambda_\ell=0$ then clearly it is a partition counted by $A(n)$.
\end{itemize}
By union of cases $1$ and case $2$, we get all partitions counted by $D_{k-1}(n-k+1)$ and $2A(n)$. This completes the proof of the theorem.

\section{Analytic Proofs}\label{analytic}

\begin{proof}[\textbf{Theorem \textup{\ref{B_k=C_k}}}][]
For $0\leq m\leq k-1$, let $C_{k,m}(n)$ denote the number of partitions of $n$ which have an even part (say $2\ell$) and parts which are not exceeding than $\ell$ are distinct, and only $m$ parts are larger than $2\ell$ which are distinct and even lying in the interval $[2\ell+2,2\ell+2k-2]$. Moreover, let $B_{k,m}(n)$ denote the number of partitions of $n$ having odd parts, and only except $m$ parts (even and distinct) which lies in the interval $[2\ell+2,2\ell+2k-2]$ ($2\ell-1$ is the largest odd part). Consider the generating function, (where $x$ counts the number of parts after $2\ell$)
\begin{align*}
\sum_{n=1}^{\infty}\frac{(-q;q)_n(-xq^{2n+2};q^2)_{k-1}}{(q^{n+1};q)_n}q^{2n}=\sum_{n=1}^{\infty}\left(\sum_{m=0}^{k-1}C_{k,m}(n)x^m\right)q^n.
\end{align*}
Now the left-hand side of the above equation can be rewritten as
\begin{align*}
\sum_{n=1}^{\infty}\frac{(-q;q)_n(-xq^{2n+2};q^2)_{k-1}}{(q^{n+1};q)_n}q^{2n}&=\sum_{n=1}^{\infty}\frac{(-xq^{2n+2};q^2)_{k-1}}{(q;q^2)_n}q^{2n}\\
&=q\sum_{n=1}^{\infty}\frac{(-xq^{2n+2};q^2)_{k-1}}{(q;q^2)_n}q^{2n-1}\\
&=q\sum_{n=1}^{\infty}\left(\sum_{m=0}^{k-1}B_{k,m}(n)x^m\right)q^n\\
&=\sum_{n=1}^{\infty}\left(\sum_{m=0}^{k-1}B_{k,m}(n)x^m\right)q^{n+1}.
\end{align*}
On comparing the coefficients of $x^{2j}q^n$ both  sides, we deduce that
\begin{align*}
B_{k,2j}(n)=C_{k,2j}(n+1).
\end{align*}
Now summing over over $j$, we obtain 
$$B_{k}^e(n)=C_{k}^e(n+1).$$
In the same way, we can prove that
		$$B_{k}^o(n)=C_{k}^o(n+1).$$

\end{proof}

\begin{proof}[\textbf{Theorem \textup{\ref{BkCkDk}}}][]
The first equality is simple consequence of Theorem \ref{B_k=C_k}.

We next prove that
\begin{align}
\frac{1}{2}D_{2k}(n)=C_k^e(n)-C_k^o(n),\quad n>2^{k-1}k(2k-1),\nonumber
\end{align}
%and
%\begin{align}\label{bkck}
%B_{k}^e(n)-B_{k}^o(n)=C_{k}^e(n+1)-C_{k}^o(n+1), \quad n>1.
%\end{align}
To prove this, observe that
\begin{align}\label{before euler}
&\frac{1}{2}\sum_{m=0}^{\infty}q^{2mk}(-q^{m+1};q)_\infty+\frac{1}{2}(q;q)_{2k-1}\nonumber\\
&=\frac{1}{2}(-q;q)_\infty\sum_{m=0}^{\infty}\frac{q^{2mk}}{(-q;q)_m}+\frac{1}{2}(q;q)_\infty\sum_{m=0}^{\infty}\frac{q^{2mk}}{(q;q)_m},
\end{align}
where we used the Euler's identity \cite[p.~19, (2.2.5)]{andrews}
\begin{align}\label{euler}
\frac{1}{(t;q)_\infty}=\sum_{m=0}^\infty\frac{t^m}{(q;q)_m},
\end{align}
with replacing $t$ by $q^{2k}$. Simplifying the right-hand side of \eqref{before euler}, we arrive at
\begin{align}\label{after euler}
&\frac{1}{2}\sum_{m=0}^{\infty}q^{2mk}(-q^{m+1};q)_\infty+\frac{1}{2}(q;q)_{2k-1}\nonumber\\
&=\frac{1}{2}(q^2;q^2)_\infty\ \sum_{m=0}^{\infty}\frac{q^{2mk}}{(q^2;q^2)_m}\left(\frac{1}{(q^{m+1};q)_\infty}+\frac{1}{(-q^{m+1};q)_\infty}\right)\nonumber\\
&=\frac{1}{2}(q^2;q^2)_\infty\ \sum_{m=0}^{\infty}\frac{q^{2mk}}{(q^2;q^2)_m}\left(\sum_{n=0}^{\infty}\frac{q^{n(m+1)}}{(q;q)_n}+\sum_{n=0}^{\infty}\frac{(-1)^nq^{n(m+1)}}{(q;q)_n}\right),
\end{align}
where the final step is obtained by applying \eqref{euler} twice, first with $t=q^{m+1}$ and then with $t=-q^{m+1}$. Equation \eqref{after euler} can be rewritten as
\begin{align*}
\frac{1}{2}\sum_{m=0}^{\infty}q^{2mk}(-q^{m+1};q)_\infty+\frac{1}{2}(q;q)_{2k-1}
&=(q^2;q^2)_\infty\sum_{m,n=0}^{\infty}\frac{1}{2}\left((1+(-1)^n)\right)\frac{q^{2km+n+nm}}{(q;q)_n(q^2;q^2)_m}\\
&=(q^2;q^2)_\infty\sum_{m,n=0}^{\infty}\frac{1}{2}\frac{q^{2km+2n+2nm}}{(q;q)_{2n}(q^2;q^2)_m}\\
&=(q^2;q^2)_\infty\sum_{n=0}^{\infty}\frac{q^{2n}}{(q;q)_{2n}}\sum_{m=0}^{\infty}\frac{q^{2nm+2km}}{(q^2;q^2)_{m}}.
\end{align*}
Employing \eqref{euler} again for the inner sum on the right-hand of the above equation, we deduce that
\begin{align*}
\frac{1}{2}\sum_{m=0}^{\infty}q^{2mk}(-q^{m+1};q)_\infty+\frac{1}{2}(q;q)_{2k-1}&=(q^2;q^2)_\infty\sum_{n=0}^{\infty}\frac{q^{2n}}{(q;q)_{2n}(q^{2n+2k};q^2)_\infty}\\
&=\sum_{n=0}^{\infty}\frac{(-q;q)_n(q^{2n+2};q^2)_{k-1}}{(q^{n+1};q)_n}q^{2n}\\
&=\sum_{n=1}^{\infty}\frac{(-q;q)_n(q^{2n+2};q^2)_{k-1}}{(q^{n+1};q)_n}q^{2n}+(q^2;q^2)_{k-1}.
\end{align*}
Now, this clearly implies that
\begin{align*}
\frac{1}{2}\sum_{n=0}^{\infty} D_{2k}(n)q^n+\frac{1}{2}(q;q)_{2k-1}=\sum_{n=0}^\infty \left(C_k^e(n)-C_k^o(n)\right)q^n+(q^2;q^2)_{k-1},
\end{align*}
Finally, comparing coefficients of $q^n$ for $n>2^{k-1}k(2k-1)$, we completes the proof of the second equality of \eqref{BkCkDk eqn}.

%Now observer that
%\begin{align*}
%\sum_{n=1}^\infty \left(C_k^e(n)-C_k^o(n)\right)q^n&=\sum_{n=1}^{\infty}\frac{(-q;q)_n(q^{2n+2};q^2)_{k-1}}{(q^{n+1};q)_n}q^{2n}\\
%&=\sum_{n=1}^{\infty}\frac{(q^{2n+2};q^2)_{k-1}}{(q;q^2)_n}q^{2n}\\
%&=\sum_{n=1}^{\infty}\left(B_k^e(n)-B_k^o(n)\right)q^{n+1}.
%\end{align*}
%This proves \eqref{bkck}. Hence, we complete the proof of the theorem.
\end{proof}

\begin{proof}[\textbf{Theorem \textup{\ref{AkDk}}}][]
Note that
\begin{align}\label{enum}
&\sum_{n=1}^{\infty}\left(P_1(d,n)+P_{k-1}'(d,n)+P_2(d,n)+P_{k-1}^{''}(d,n)\right)q^{n+1} \nonumber\\
&=q\sum_{n=1}^{\infty}P_1(d,n)q^{n}+q\sum_{n=1}^{\infty}P_{k-1}
'(d,n)q^{n}+q\sum_{n=1}^{\infty}P_2(d,n)q^{n}+q\sum_{n=1}^{\infty}P_{k-1}^{''}(d,n)q^{n} \nonumber\\
&=q\sum_{n=1}^{\infty}q^{n+1}(-q^{n+2};q)_\infty+q^{k}(-q^{2};q)_\infty+q\sum_{n=1}^{\infty}q^n(-q^{n+2};q)_\infty+q\sum_{n=1}^{\infty}q^{n+(n+1)(k-1)}(-q^{n+2};q)_\infty \nonumber\\
&=\sum_{n=1}^{\infty}(1+q)q^{n+1}(-q^{n+2};q)_\infty+q^{k}(-q^{2};q)_\infty+\sum_{n=2}^{\infty}q^{nk}(-q^{n+1};q)_\infty \nonumber\\
&=\sum_{n=1}^{\infty}(1+q)q^{n+1}(-q^{n+2};q)_\infty+\sum_{n=1}^{\infty}q^{nk}(-q^{n+1};q)_\infty.
%&=\sum_{n=1}^{\infty}D_{k}(n)q^n.
\end{align}
Now observing that the first term in the ultimate line is the generating function for the number of partitions of $n>1$ into distinct parts, whereas the second term enumerates the partitions of $n$ in which the smallest part is repeated $k$ times and all other parts are distinct. Together, these two classes account for the partitions enumerated by $D_k(n)$.

Now upon comparing the coefficients on both sides of \eqref{enum}, we complete the proof of our claim in \eqref{AkDk eqn}.
\end{proof}

\begin{proof}[\textbf{Corollary \textup{\ref{special case k=1}}}][]
Let $k=1$ in \eqref{genral result eqn}. Using the facts that
$B_1^e(n)=B(n),\ B_1^o(n)=0=C_1^o(n)$ and $C_1^e(n)=C(n)$ and $A_2(n)=A(n)$, $D_2(n)=D(n)$, we arrive at \eqref{main theorem eqn}.

\end{proof}

\begin{proof}[\textbf{Theorem \textup{\ref{extd euler}}}][]
As it is already known, from Theorem \ref{main theorem}, that $ A(n)=B(n)=C(n+1)=\frac{1}{2}D(n+1)$, it is enough to show that $C(n)=F(n)$ and  $C(n)=E(n+1)$. To that end, note that
\begin{align*}
\sum_{n=1}^{\infty}C(n)q^n&=\sum_{n=1}^{\infty}\frac{(-q;q)_n}{(q^{n+1};q)_n}q^{2n}\\
&=\sum_{n=1}^{\infty}\frac{(q^2;q^2)_n}{(q;q)_{2n}}q^{2n}\\
&=\sum_{n=1}^{\infty}\frac{q^{2n}}{(q;q^2)_n}\\
&=\sum_{n=1}^{\infty}F(n)q^n.
\end{align*}

Next consider,
\begin{align*}
\sum_{n=0}^{\infty}E(n)q^n	&=\sum_{n=0}^{\infty}\frac{q^{2n+1}}{(q;q^2)_n}\\
&=q\sum_{n=0}^{\infty}\frac{q^{2n}}{(q;q^2)_n}\\
&=q+q\sum_{n=1}^{\infty}\frac{q^{2n}}{(q;q^2)_n}\\
&=q+\sum_{n=1}^{\infty}C(n)q^{n+1}.
\end{align*}

\noindent
\textbf{Bijective proof of $B(n)=F(n+1)$ and $B(n)=E(n+2)$:}		
We obtain a simple bijection between $B(n)$ and $F(n+1)$ by adding 
$1$ to the largest part of a partition counted by $B(n)$ and, conversely, subtracting $1$ from the largest part of a partition counted by $F(n+1)$. Similarly, a bijection between $B(n)$ and $E(n+2)$ is obtained by adding $2$ to, and subtracting $2$ from, the largest part of the corresponding partitions.
\end{proof}

\begin{proof}[\textbf{Theorem \textup{\ref{legendre analogue}}}][]
It is clear that 
\begin{align}\label{k-1}
\sum_{n=1}^\infty\left(D_k^e(n)-D_k^o(n)\right)q^n&=\sum_{n=0}^\infty q^{nk}\left(q^{n+1};q\right)_\infty\nonumber\\
&=(q;q)_{k-1},
\end{align}
where we employed the $q$-binomial theorem \cite[p.~17, Theorem 2.1]{andrews}
\begin{align*}
\sum_{n=0}^\infty\frac{(a;q)_n}{(q;q)_n}z^n=\frac{(az;q)_\infty}{(z;q)_\infty}.
\end{align*}
Let $n\leq k-1$. Observe that the coefficient of $q^n$ on the right-hand side of \eqref{k-1} is $P_e(d,n)-P_o(d,n)$. Hence, on comparing coefficients of $q^n$ on both sides of \eqref{k-1}, we have
\begin{align}\label{equal}
D_k^e(n)-D_k^o(n)=P_e(d,n)-P_o(d,n), \ n\leq k-1.
\end{align}
Now, for $1\leq n\leq k-1$ of the form $n=\frac{1}{2}m(3m\pm1)$, applying \eqref{legendre} in \eqref{equal}, we prove the first part. 

Now let $k-1< n\leq \frac{1}{2}k(k-1)$. Observe that coefficient of $q^n$ on the right-hand side of \eqref{k-1} is $P_e^{(k)}(d,n)-P_o^{(k)}(d,n)$. Hence, comparing the coefficients of \eqref{k-1}, we arrive at the second part. 

The last part is easy as the coefficients of $q^n$ on the right-hand side of \eqref{k-1} are zero for every $n>\frac{1}{2}k(k-1)$. 
This completes the proof.
\end{proof}

\begin{proof}[\textbf{Corollary \textup{\ref{parity of D_k}}}][]
Let $n>\frac{1}{2}k(k-1)$. The result in Corollary is quite straightforward from the third part of \eqref{legendre analogue eqn}.
\end{proof}

\begin{proof}[\textbf{Theorem \textup{\ref{finite analogue}}}][]
Note that \eqref{finite analogue eqn} can be rephrased as
\begin{align*}
\sum_{n=0}^{N}\frac{q^{kn}}{(-q;q)_n}&=\sum_{i=0}^{k-1}(-1)^{k-1+i}(q^{i+1};q)_{k-1-i}\left(2-\frac{q^{(N+1)i}}{(-q;q)_N}\right).
\end{align*}

We will prove this by induction on $k$ for any non-negative integer $N$. Note that our result in \eqref{finite analogue eqn}  holds for $k=1$  due to the result of Andrews given in \eqref{andrews amm}. Now suppose \eqref{finite analogue eqn} is true for $k-1$. Next, consider
\begin{align*}
\sum_{n=0}^{N}\frac{q^{kn}}{(-q;q)_n}&=	\sum_{n=0}^{N}\frac{q^{(k-1)n}(-1+(1+q^n))}{(-q;q)_n}\\
&=-\sum_{n=0}^{N}\frac{q^{(k-1)n}}{(-q;q)_n}+	\sum_{n=0}^{N}\frac{q^{(k-1)n}(1+q^n)}{(-q;q)_n}\\
&=-\sum_{n=0}^{N}\frac{q^{(k-1)n}}{(-q;q)_n}+\sum_{n=1}^{N}\frac{q^{(k-1)n}(1+q^n)}{(-q;q)_n}+2\\
&=-\sum_{n=0}^{N}\frac{q^{(k-1)n}}{(-q;q)_n}+\sum_{n=1}^{N}\frac{q^{(k-1)n}}{(-q;q)_{n-1}}+2.
\end{align*}
Replacing $n$ by $n+1$ in the second sum on the right-hand side and then adding and subtracting the corresponding term with $n=N+1$, we obtain
\begin{align*}
\sum_{n=0}^{N}\frac{q^{kn}}{(-q;q)_n}&=-\sum_{n=0}^{N}\frac{q^{(k-1)n}}{(-q;q)_n}+q^{k-1}\sum_{n=0}^{N}\frac{q^{(k-1)n}}{(-q;q)_n}-\frac{q^{(k-1)(N+1)}}{(-q;q)_N}+2\\
&=-(1-q^{k-1})\sum_{n=0}^{N}\frac{q^{(k-1)n}}{(-q;q)_n}-\frac{q^{(k-1)(N+1)}}{(-q;q)_N}+2.
\end{align*}
Now, applying \eqref{finite analogue eqn} for $k-1$, which is valid by the induction hypothesis, we obtain
\begin{align*}
\sum_{n=0}^{N}\frac{q^{kn}}{(-q;q)_n}&=-(1-q^{k-1})\left(\sum_{i=0}^{k-2}(-1)^{k-2+i}(q^{i+1};q)_{k-2-i}\left(2-\frac{q^{(N+1)i}}{(-q;q)_N}\right)\right)-\frac{q^{(k-1)(N+1)}}{(-q;q)_N}+2\\
&=\sum_{i=0}^{k-2}(-1)^{k-1+i}(q^{i+1};q)_{k-1-i}\left(2-\frac{q^{(N+1)i}}{(-q;q)_N}\right)+\left(2-\frac{q^{(k-1)(N+1)}}{(-q;q)_N}\right)\\
&=\sum_{i=0}^{k-1}(-1)^{k-1+i}(q^{i+1};q)_{k-1-i}\left(2-\frac{q^{(N+1)i}}{(-q;q)_N}\right).
\end{align*}
This proves \eqref{finite analogue} for general $k$, hence we complete the proof.
\end{proof}

\begin{proof}[\textbf{Theorem \textup{\ref{D_k in terms of polyn}}}][]
 Note that \eqref{finite analogue eqn} can be rewritten as
 \begin{align*}
&\sum_{n=0}^{N}q^{kn}(-q^{n+1};q)_\infty\\
&=2(-q;q)_\infty\sum_{j=0}^{k-1}(-1)^{k-1+j}(q^{j+1};q)_{k-1-j}-\frac{(-q;q)_\infty}{(-q;q)_N}\sum_{j=0}^{k-1}(-1)^{k-1+j}(q^{j+1};q)_{k-1-j}q^{(N+1)j}\\
 &=2(-q;q)_\infty\sum_{j=0}^{k-1}(-1)^{k-1+j}(q^{j+1};q)_{k-1-j}-(-1)^{k-1}(q;q)_{k-1}\frac{(-q;q)_\infty}{(-q;q)_N}\\
 &\qquad-(-q;q)_\infty\sum_{j=1}^{k-1}(-1)^{k-1+j}(q^{j+1};q)_{k-1-j}\frac{q^{(N+1)j}}{(-q;q)_N}.
 \end{align*}
Taking $N\to\infty$ in the above equation and noting that the last finite sum vanishes, since $\lim_{N\to\infty}q^{(N+1)j}=0$ for $j\geq1$, we obtain
\begin{align*}
&\sum_{n=0}^{\infty}q^{kn}(-q^{n+1};q)_\infty\\
&=2(-q;q)_\infty\sum_{j=0}^{k-1}(-1)^{k-1+j}(q^{j+1};q)_{k-1-j}+(-1)^{k}(q;q)_{k-1}\\
 &=2(-q;q)_\infty\sum_{j=0}^{k-1}(-1)^{k-1+j}(q^{j+1};q)_{k-1-j}+(-1)^{k}(q;q)_{k-1}.
 \end{align*}
Finally, the result follows upon replacing $j$ by $k-1-j$ in the finite sum on the right-hand side of the above expression.
 \end{proof}

\begin{proof}[\textbf{Theorem \textup{\ref{D_k+D_k}}}][]
Observe that \eqref{D_k in terms of polyn eqn} can be rewritten as
\begin{align*}
\sum_{n\ge0}^{}{D_k(n)q^n}
&=2(-q;q)_{\infty}+2(-q;q)_\infty(1-q^{k-1})\sum_{j=1}^{k-1}(-1)^j(q^{k-j};q)_{j-1}+(-1)^k(q;q)_{k-1}\\
&=2(-q;q)_{\infty}-2(-q;q)_\infty(1-q^{k-1})\sum_{j=0}^{k-2}(-1)^j(q^{k-j-1};q)_{j}-(-1)^{k-1}(1-q^{k-1})(q;q)_{k-2}\\
&=2(-q;q)_{\infty}-\left\{2(-q;q)_\infty\sum_{j=0}^{k-2}(-1)^j(q^{k-j-1};q)_{j}+(-1)^{k-1}(q;q)_{k-2}\right\}\\
&\qquad+q^{k-1}\left\{2(-q;q)_\infty \sum_{j=0}^{k-2}(-1)^j(q^{k-j-1};q)_{j}+(-1)^{k-1}(q;q)_{k-2}\right\}\\
&=2\sum_{n\ge0}^{}A(n))q^n-\sum_{n\ge0}^{}D_{k-1}(n)q^n+\sum_{n\ge0}^{}D_{k-1}(n))q^{n+k-1},
\end{align*}
where the last step follows upon employing \eqref{D_k in terms of polyn eqn} twice. Now, comparing coefficients of $q^n$ on both sides completes the proof.
\end{proof}

%%\begin{proof}[\textbf{Corollary \textup{\ref{A_n=D_n}}}][]
%%Letting $k=2$ in \eqref{D_k in terms of polyn eqn} and using the fact that $D_2(n)=D(n)$, we see that
%%\begin{align*}
%%\sum_{n=1}^\infty D(n)q^n=2q(-q;q)_\infty+(1-q).
%%\end{align*}
%%Now, comparing the coefficient of $q^n$ for $n>1$, we prove the corollary.
%%\end{proof}

\begin{proof}[\textbf{Corollary \textup{\ref{D3A}}}][]
Substituting $k=3$ in \eqref{D_k in terms of polyn eqn}, we obtain
\begin{align*}
\sum_{n=0}^{\infty}D_3(n)q^n&=((q^{2}-1)2q+2)(-q;q)_\infty-(1-q)(1-q^2)\\
&=(2q^3-2q+2)(-q;q)_\infty-1+q^2+q-q^3\\
&=2\sum_{n=3}^{\infty}A(n-3)q^n-2\sum_{n=1}^{\infty}A(n-1)q^n+2\sum_{n=0}^{\infty}A(n)q^n-1+q^2+q-q^3.
\end{align*}
Equation \eqref{D3A eqn} now follows upon comparing coefficients for $n\geq4$.
\end{proof}

\section{Concluding Remarks}\label{concluding remarks}
In this work we have introduced and studied several classes of partitions, including those involving the parity of parts, differences of partition numbers, and partitions with a repeated smallest part. These types of partition classes have been extensively investigated in the literature. We have provided analytic proofs of all our results, together with bijective proofs for many of them.

We propose the following questions for further investigation.
\begin{enumerate}
\item We provided bijection for every equality in \eqref{genral result eqn} except for those involving $B_k^e(n)-B_k^o(n)$ or $C_k^e(n+1)-C_k^o(n+1)$ and $\frac{1}{2}D_{2k}(n+1)$ or $A_{2k}(n)$. It would therefore be interesting to construct explicit bijections between these classes of partitions.

\item By decomposing $D(n)$ into two classes $D_2^e(n)$ and $D_2^o(n)$, we reformulated Theorem \ref{main theorem} as given in  \eqref{main theorem re}. It would be great to find bijections between the sets counted by $D_2^e(n+1)$ (or $D_2^o(n+1)$) and those counted by  $A(n), B(n)$ or $C(n)$. The same question applies for general $k$ in \eqref{genral result re}.

\item Theorem \ref{legendre analogue} gives the following identity
\begin{align*}
D_k^e(n)-D_k^o(n)= P_e^{(k)}(d,n)-P_o^{(k)}(d,n), \quad \textup{if}\ k-1<n\leq\frac{1}{2}k(k-1),\ k>2.
\end{align*}
Can one find a direct bijective proof of this equality?

\end{enumerate}

\medskip

\noindent {\bf{Acknowledgements:}}\,
% We are grateful to the referees for their numerous helpful comments and suggestions which have greatly improved the exposition of this paper. We also thank Ae Ja Yee and Atul Dixit for helpful discussions. 
Authors would like to thank Atul Dixit and Aritram Dhar for their valuable comments on an earlier version of the paper. 
The first author is partially supported by Grant ANRF/ECRG/ 2024/003222/PMS of the Anusandhan National Research Foundation, Govt. of India, and  the FIG grants of IIT Roorkee. The second author is supported by UGC fellowship for her PhD. Both the authors sincerely thank these institutions for their support.

\end{document}